\documentclass[12pt]{article}

\usepackage{amssymb,latexsym}
\usepackage{authblk}

\newcommand{\ZZ}{\mathbb{Z}}
\newcommand{\NN}{\mathbb{N}}
\newcommand{\RR}{\mathbb{R}}

\newtheorem{thm}{Theorem}
\newtheorem{defi}{Definition}
\newtheorem{lem}{Lemma}

\newtheorem{cor}{Corollary}
\newtheorem{pro}{Proposition}
\newtheorem{rem}{Remark}
\newtheorem{ques}{Questions}

\title{Equalities, inequalities between density flow, spatial and temporal entropies of cellular automata}

\author[1]{Pierre Tisseur }
\affil[1]{Trent University 

Department of Mathematics, Peter Gzowsky College, Peterborough, Ontario, Canada}
\affil[1]{ Centro de Matematica, Computa\c{c}\~ao e Cogni\c{c}\~ao, Universidade Federal do ABC, Santo Andr\'e, S\~ao Paulo, Brasil \thanks{E-mail address: \texttt{pierre.tisseur@ufabc.edu.br}}}

\date{ June 2006}

\begin{document}

\maketitle

\begin{abstract}
Defining the density flow of perturbations moving at a given speed for cellular automata, we establish
equalities and inequalities between the measurable entropy of a cellular automaton and
  the measurable entropy of its associated shift.
  We illustrate our results by different examples and we study some relations
  between the density flow (with respect to a  linear or  a sublinear speed) and some properties of
  cellular automata
 ( positive expansiveness, $\mu$-expansiveness and  $\mu$-equicontinuity).
The probability measure we consider must be shift ergodic and invariant for the automaton.
\end{abstract}

\section{Introduction}
A one-dimensional cellular automaton (CA) is a discrete mathematical idealization of a space-time
physical system. The space, called configuration space, is the set of
doubly infinite sequences of elements of a finite set $A$. The discrete time is
represented by the action of a cellular automaton $F$ on this space.
Since  cellular automata are  continuous, shift commuting dynamical systems and depend on a local rule,
 relations between the spacial (shift) and temporal (automaton)
entropies seems natural. In \cite{Sh92} Shereshevsky establishes a first relation between the entropy of
the shift and the entropy of the cellular automaton using some discrete analog of Lyapunov
exponents for a CA ergodic measure. The product of the shift entropy by the
lyapunov exponents give an upper bound of
the CA.
 In \cite{ti2000} Tisseur define average Lyapunov exponents  using a CA invariant and shift ergodic
measure.
Given an infinite configuration, the two Lyapunov exponents (left and right) represent the speed of the
faster perturbations
moving from the left to the right or from the right to the left coordinates.
Remark that it is possible to define these left and right exponents for each points or to take an
average value of these speeds.

 The Average Lyapunov exponents give the rate of propagation of the faster perturbations  and do not take in
account the amount of information the automaton carry. This is  the reason for which in many cases the
inequality is not an equality.
Nevertheless  the  average Lyapunov exponents (when they are equal to zero) are able to prove that the
measurable entropy of the cellular automaton  is null for CA  with equicontinuous points
(see \cite{ti2000}) and for some particular sensitive (without equicontinous point) CA (see \cite{BT2006}).

Here, we propose to introduce the density flow $M_\mu(v)$ of the  perturbations
with respect to  a velocity $v=(v^+,v^-)$ where $v^+$ is the right to left speed
and $v^-$ the right to left speed.
This average flow  depends on a  shift-ergodic and CA invariant probability measure $\mu$.
In section 3.2 (Theorem 2), we establish that
$h_\mu (F)=h_\mu (\sigma )\times M_\mu (v)\times (v^++v^-)$ where $h_\mu (F)$ and
$h_\mu (\sigma )$ are respectively the entropy of the automaton $F$ and the entropy of the
shift $\sigma$ and the speeds $v^+$ and  $v^-$ are greater than
the radius $r$ of the (CA). More generally, Theorem 2  states that  $h_\mu (F)\ge h_\mu (\sigma )\times M_\mu (v)\times (v^++v^-)$.
These relations can be compared with the inequality that appears in \cite{ti2000} :
$h_\mu (F)\le h_\mu (\sigma )\times (I_\mu^++I_\mu^-)$ where $I_\mu^+$ and $I_\mu^-$
are the left and right Average Lyapunov exponents.
Remark that the condition $v^+,v^-\ge r$ can be substituted by a weaker one (see Theorem 2)
in relation with "maximum Lyapunov exponents".

One of our  motivations to define the density flow is to understand the dynamic of cellular
automata with null measurable entropy.
From \cite{ti2000}, $I_\mu^++I_\mu^-=0$ and $h_\mu (F)=0$, if there exist equicontinuous points
in the topological support of the shift ergodic and $F$-invariant measure $\mu$.
Having equicontinuous points is equivalent to say that exist some  patterns of positive measure, called blocking
 words that stop the perturbations.
 In \cite{BT2006}, Bressaud and Tisseur show that there exists a sensitive cellular automaton (without equicontinuous points)
 such that $I_\mu^++I_\mu^-=0$ which implies  that $h_\mu (F)=0$.
 In section 3.3 we give an example of a cellular automaton $F$,  such that the faster perturbations move at
 positive speed ($I_\mu^++I_\mu^-=1$) but the density flow at speed $v^+,v^-:=1$ is equal to zero.
 In that case the density flow is equal to zero because the weight of the perturbation moving at
 a positive speed is equal to zero and
 Theorem 2 implies that $h_\mu (F)=0$.

Finally  we show some relations between some subclasses of CA and the density of flow.
We state that there always exists a positive velocity $v$ such that $M_\mu (v)=1$ when
the cellular automaton $F$ has the positive expansiveness property  ($\mu$ is
a shift ergodic and $F$-invariant measure).
In this case,  all the perturbations move with a speed at least equal to $v$.
Furthemore, we prove that $M_\mu (v)=0$ for all linear or "sublinear speed" $v$ if
the cellular automaton $F$ has  $\mu$-equicontinuous points, a measurable equivalent
to the existence of equicontinuous points introduced by Gilman in \cite{GI87}.
This class is the complementary class of the $\mu$-expansive class which is a
measurable equivalent to the positive expansiveness one.
 We call sublinear speed $v$ any couple $(v^+,v^-)$ where
$v^+$ and $v^-$ represent  positive integer sequences $(v_n^+)$ and $(v^-_n)$ that
verify $\lim_{n\to\infty}(v_n^++v_n^-)=+\infty$ and $\lim_{n\to\infty}(\frac{v_n^+}{n}+\frac{v_n^-}{n})=0$.



\section{Preliminary}

\subsection{Symbolics systems and  cellular automata}

Let $A$ be a finite set or alphabet. Denote by $A^*$ the set of all
concatenations of letters in $A$. These concatenations are called words. The
length of a word $u\in A^*$ is denoted by  $\vert u\vert$.
The set of bi-infinite sequences
$x=(x_i)_{i\in\ZZ}$ is denoted by $A^\ZZ$. A point $x\in A^\ZZ$ is called a
configuration. For $i\le j$ in $\ZZ$ we denote by $x(i,j)$ the word $x_i\ldots
x_j$ and by $x(p,\infty )$ the infinite sequence $(v_i)_{i\in\NN}$ such that for
all $i\in\NN$ one has $v_i=x_{p+i-1}$. We endow $A^\ZZ$ with the product
topology. The shift $\sigma \colon A^\ZZ\to A^\ZZ$ is defined by :
$\sigma (x)=(x_{i+1})_{i\in \ZZ}$. For each integer $t$ and each word $u$, we
call cylinder the set $[u]_t=\{x\in A^\ZZ : x_t=u_1\ldots ;x_{t+\vert
u\vert}=u_{\vert u\vert}\}$. For this topology $A^\ZZ$ is a compact metric
space. A metric compatible with this topology can be defined by the distance
$d(x,y)=2^{-i}$ where $i=\min\{\vert j\vert \,\mbox{ such that } x(j)\ne y(j)\}$.
The dynamical system $(A^\ZZ ,\sigma )$ is called the full shift. A subshift $X$
is a closed shift-invariant subset $X$ of $A^\ZZ$ endowed with the shift
$\sigma$. It is possible to identify $(X,\sigma )$ with the set $X$.
 If $\alpha =\{A_1,\ldots ,\,
A_n\}$ and $\beta =\{B_1,\ldots ,\, B_m\}$ are two partitions denote
by $\alpha \vee \beta$ the partition $\{A_i\cap A_j\, i=1\ldots n, \,\,
j=1,\ldots ,\, m\}$.
Consider a probability measure $\mu$ on the Borel sigma-algebra
$\cal{B}$ of $A^\ZZ$. If $\mu$ is $\sigma$-invariant then the topological
support of $\mu$ is a subshift denoted by $S(\mu )$.
The metric entropy $h_\mu (T)$ of a transformation $T$
is an isomorphism invariant between two $\mu$-preserving
transformations.
Put $H_\mu (\alpha ) = \sum_{A\in\alpha}\mu (A)\log \mu (A)$.
The entropy of the partition $\alpha$ is defined as $h_\mu (\alpha ) =
\lim_{n\to\infty}1/nH_\mu (\vee_{i=0}^{n-1}T^{-i}\alpha )$ and the entropy of $(X,T,\mu )$ as
$\sup_\alpha h_\mu (\alpha )$.

 A cellular automaton  (CA) is a continuous self-map $F$ on
$A^\ZZ$ commuting with the shift. The Curtis-Hedlund-Lyndon theorem
states  that for every cellular automaton $F$ there exist an integer $r$ and a
block map $f$ from $A^{2r+1}$ to $A$ such that: $F(x)_i=f(x_{i-r},\ldots ,x_i
,\ldots ,x_{i+r}).$
 The integer $r$ is called the radius of the cellular
automaton. If the block map of a cellular automaton is such that
$F(x)_i=f(x_{i},\ldots ,\ldots ,x_{i+r})$,
 the cellular automaton is called
one-sided and can be extended a map on a two-sided shift $A^\ZZ$ or a map on a
one-sided shift $A^\NN$. If $X$ is a subshift of  $A^\ZZ$ and one has
$F(X)\subset X$, the restriction of $F$ to $X$ determines a dynamical system
$(X,F)$; it is
called a cellular automaton  on $X$.

\subsection{Other definitions and one inequality}

\paragraph{Topological and measurable properties:}
$\mbox{  }$
\medskip

For all $x\in A^\ZZ$ and $\epsilon >0$ denote by $D(x,\epsilon )$ be the set of all points $y$ such that
 for all $i\in \NN$ one has $d(F^i(x),F^i(y))<\epsilon$ and for all positive integer $n$, write
 $B_n(x)=D(x,2^{-n})$. For all $x\in A^\ZZ$, the set $C_n(x)$ represents the set of all points
 $y$ such that $y(-n,n)=x(-n,n)$.

\begin{defi} Equicontinuity
$\mbox{ }$

- A point $x$ is called an equicontinuous point if for all positive integer $n$, there exists another
integer $m\ge n$ such that $B_n(x)\supset C_n(x)$.

- A point $x$ is called a $\mu$-almost equicontinuous point if $\mu (B_n(x))>0$ for all integer
$n\ge r$.

\end{defi}

From \cite{GI87}, and \cite{tippt2006}, if $F$ is a CA of radius $r$,  $\mu$ is a shift ergodic measure
and  $x$ a point which verifies  $\mu (B_r(x))>0$  then there exist a set of measure one of
$\mu$-equicontinuous
points.
In this case we say that $F$ is a {\it $\mu$-equicontinuous} CA.

\begin{defi} Expansiveness
$\mbox{ }$

-A Cellular automaton is positively expansive if there exists a positive integer  $n$
such that for all $x\in A^\ZZ$ one has $B_n(x)=\{x\}$.

-Let $\mu$ be a shift ergodic measure. A  cellular automaton  $F$ is $\mu$-almost expansive if there exists
a positive integer  $n$ such that for all $x\in A^\ZZ$, $\mu (B_n(x))$ $=0$ (see \cite{GI87}, \cite{tippt2006}).

-A cellular automaton is sensitive if for all points $x\in A^\ZZ$ and each
 integer $m>0$ one has $C_m(x)\subsetneq B_r(x)$.
\end{defi}

\paragraph{Lyapunov exponents :}
$\mbox{ }$
\medskip



The Average Lyapunov exponents represent the average speed of the faster perturbations on  the
one dimensional lattice $A^\ZZ$.

Let $W_i^+(x)$ and $W_i^-(x)$ are respectively the set of all the  perturbations
(any element in $W_i^\pm (x)$ different of $x$)
 at the left side and at the right side  of a coordinate $i$ in a point $x$.
More formaly
$$
W_s^+(x)=\{y\in A^\ZZ :\,\forall i\ge s, \, y_i=x_i \},\hskip .2 cm
        W_s^-(x)=\{y\in A^\ZZ :\,\forall i\le s, \, y_i=x_i\}.
$$

Let's define two continuous sequences of functions $I^-_n$ and $I^+_n$ from
$A^\ZZ\to \NN$ by

\[
I^-_n(x)=\min\{ s\in\NN \, \vert\, \, \forall \, 1\le i\le n ,\,\vert \,
F^i(W^-_{s}(x))\subset W_0^-(F^i(x))\},
\]
\[
I^+_n(x)=\min\{ s\in\NN \,\vert \, \, \forall \, 1\le i\le n ,\,\vert\,
F^i(W^+_{-s}(x))\subset W_0^+(F^i(x))\}.
\]

Remark that
$I^+_n$ and $I^-_n$ are  bounded by $rn$ where $r$ is the radius of the CA.
Set $
I^+_{n;\mu}=\int_X I^+_n(x)d\mu (x)$ and $I^-_{n;\mu}=
\int_X I^-_n(x)d\mu(x).$ The Birkhoff's theorem implies that for almost all $x$ one has
 $I^+_{n;\mu}=\lim_{n\to\infty}\sum_{i=-m}^m\frac{1}{2m+1}I^+_n(\sigma^i
(x))$ and $I^-_{n;\mu}=\lim_{n\to\infty}\sum_{i=-m}^m\frac{1}{2m+1}I^-_n(\sigma^i(x))$.
\begin{defi}
Call average Lyapunov exponents the limits
\[
I^+_\mu=\liminf_{n\to\infty}\frac{I^+_{n;\mu}}{n}\,\, \mbox{ and }\,\,
I^-_\mu=\liminf_{n\to\infty}\frac{I^-_{n;\mu}}{n}.
\]
\end{defi}

\begin{thm}\cite{ti2000}
Let $F$ be a cellular automaton on $A^\ZZ$ and  $\mu$ is a shift ergodic and $F$-invariant  measure,
then $h_\mu (F)\le h_\mu (\sigma )(I_\mu^++I_\mu^-)$ where $h_\mu (F)$ and $h_\mu (\sigma )$ are respectively
the measurable entropies of the automaton and the shift $\sigma$.
\end{thm}


\section{Results}

This aim of this work is double:

Firstly,   establish an equality between the spatial and the temporal measurable entropy
of a cellular automaton $F$ for a shift ergodic and $F$-invariant measure.
In order to do that, we introduce the density flow with respect to a velocity $v$.

Secondly,   to  progress in the understanding  of the dynamic of cellular automata with nul measurable
entropy.

\subsection{The equality for the uniform measure}

Let's introduce the density flow of perturbations moving at  velocity $r$
in the particular case of the uniform measure.

\subsubsection{A first equality}

For each cellular automaton $F$ and positive integer $p$, denote by  $\alpha_p$
the partition of the set $A^\ZZ$ by the $p$
 central coordinates and by  $^{\{n\}}\alpha_{p}^{F}$ the partition
 $\alpha_p\vee F^{-1}\alpha_p\ldots F^{-n+1}\alpha_p$. Let  $^{\{n\}}\alpha_{p}^{F}(x)$ be the
 element of the partition $^{\{n\}}\alpha_{p}^{F}$ that contains the point $x$.
 Likewise define $^{\{\pm n\}}\alpha_p^{\sigma}(x)$ as the element of the partition
$\alpha_p\vee \sigma^{-1}\alpha_p\ldots \vee\sigma^{-n+1}\alpha_p\vee \sigma \alpha_p\vee \sigma^2\alpha_p
\ldots \vee \sigma^{n-1}\alpha_p$ generated by the shift $\sigma$ that
contains $x$.
 We suppose that $F$ is a surjective CA which implies that the
 uniform measure $\mu$ on $A^\ZZ$ is $F$-invariant (see \cite{BKM97}).
  From the probabilistic version of the Shannon-McMillan-Breiman Theorem
 we have
$h_{\mu} (F,\alpha_p)=\int_{A^\ZZ}\lim_{n\to\infty}\frac{-1}{n}\log \mu (^{\{n\}}\alpha_p^{F}(x))d\mu (x)$,
where $h_\mu (F, \alpha_p)$ is the measurable entropy of the cellular automaton $F$ with
respect to the finite partition $\alpha_p$.
 Since the cellular automaton $F$ is defined thanks to a local rule acting
on words of size $2r+1$,
 the set $^{\{n\}}\alpha_p^{F}(x))$ is a finite union of cylinders $[y(-rn-p,rn+p)]_{-rn-p}$ where
$y\in ^{\{n\}}\alpha_p^{F}(x)$. Let $T_{n,p}^{rn}(x)$ be the set of such cylinders.
Since $\mu$ is the uniform measure, all the cylinders $[x(-rn-p,rn+p)]_{-rn-p}$ have the same measure
$\vert A\vert^{-(2(p+rn)+1)}$ and
$\mu (^{\{n\}}\alpha_p^{F}(x))=\# T_{n,p}^{rn}\times \mu (^{\{rn\}}\alpha_p^{\sigma}(x))=
\# T_{n,p}^{rn}\times \vert A\vert^{-(2(p+rn)+1)}$.

Therefore, it follows that  
$$
h_{\mu} (F,\alpha_p)=\int_{A^\ZZ}\lim_{n\to\infty}\frac{-1}{n}\log \mu (^{\{n\}}\alpha_p^{F}(x))d\mu (x)
$$
$$
=\int_{A^\ZZ} \lim_{n\to\infty}\frac{-\log [\mu (^{\{\pm n\}}\alpha_p^{\sigma}(x))\times
\# T_{n,p}^{rn}(x) ]}{n}d\mu (x)
$$
$$
=\int_{A^\ZZ}\lim_{n\to\infty}\frac{\log \left(\mu (^{\{\pm rn\}}\alpha_p^{\sigma}(x))\right)}{2rn+2p+1}
\times 2r\left[1-\frac{\log (\# T_{n,p}^{rn}(x))}{-\log (\mu (^{rn}\alpha_p^{\sigma}(x)))} \right]d\mu (x).
$$

Since the uniform measure is shift ergodic, the extended version of the Shannon-McMillan-Breiman Theorem
 (see \cite{Or83}) tell us that
$$
\int_{A^\ZZ}\lim_{n\to\infty}\frac{\log \left(\mu (^{\{\pm rn\}}\alpha_p^{\sigma}(x))\right)}{2rn+2p+1}d\mu (x)=
h_\mu (\sigma, \alpha_p).
$$

This implies  that   $M^p(r):=\int_{A^\ZZ}\lim_{n\to\infty}\left[1-\frac{\log (\# T_{n,p}^{rn}(x))}
{-\log (\mu (^{\{\pm n\}}\alpha_p^{\sigma}(x)))} \right]d\mu (x)$ exits and  we obtain that
$h_{\mu}(F,\alpha_p)=h_{\mu}(\sigma, \alpha_p )\times 2r\times  M^p(r)$.
\medskip

Arguing that  $(\alpha_p)$ is a generating sequence for the shift $\sigma$,  we obtain that
$$
h_{\mu}(F)=h_{\mu}(\sigma)\times 2r\times  M^*(r),
$$
 where  $M^*(r)=\sup_{\alpha_p}M^p(r)$.
As $\mu$ is the uniform measure, $\mu (^{\{\pm n\}}\alpha_p^{\sigma}(x)))=A^{-2(rn+p)-1}$.
The terms $\frac{\log (\# T_{n,p}^{rn}(x))}{-\log (\mu (^{\pm rn}\alpha_p^{\sigma}(x)))}$ corresponds to the
ratio of the logarithm of the number of patterns $u$ of lenght $2(rn+p)+1$
 which verify that if $y\in [u]_{-rn-p}$ then $F^i(y)(-p,p)=F^i(x)(-p,p)$ for $0\le i\le n$ by
the logarithm of the number of all patterns of size $2(rn+p)+1$.
The term $M^*(r)$ can be seen as a the limit of the logarithmic proportion of the finite configurations
which move at speed $r$.

For these reasons, we call $M^*(r)$, the density flow of perturbations moving  at speed $r$.


\subsubsection{A first basic example}

Let $\sigma_1$ be the shift map on $X_1=\{0,1\}^\ZZ$ and $\sigma_2$ the $n$
iterated shift  map on $X_2=\{0,1\}^\ZZ$. For all $(x^1,x^2)\in X_1\times X_2=:X$,
one has $(\sigma_1(x^1))_i=x^1_{i+1}$ and $(\sigma_2^r (x^2))_i=x^2_{i+r}$.
Denote by $F_e$ the cellular automaton on $\{0,1\}^\ZZ\times \{0,1\}^\ZZ$ defined by
 : $F_e=\sigma_1\times\sigma_2^r$
, by $\mu$ be the uniform measure on $\{0,1\}^\ZZ\times \{0,1\}^\ZZ$ and by $\sigma$, the
shift $\sigma=\sigma_1\times \sigma_2$ on $X$.
Remark that the uniform measure $\mu$ is shift and $F$-ergodic and that  $\#^pT_n(x)$
does not depends on $x$ and $p$. Since during the $n$ first iterations of the
automaton $F$, the coordinates between $p+r+1$ and $p+rn$ in $X_1$ and the coordinates between
$-p$ and $-p-rn$ in $X$  does not reach the central coordinates between
$-p$ and $p$, we can write that  $\# T_{n,p}(x)=2^{(r-1)n+2rn}$.
We have   $\mu (\alpha_p^{rn,\sigma}(x))=4^{-2rn+2p+1}=2^{-4rn+4p+2}$ and
 $M^*(r)=\lim_{n\to\infty}1-\frac{(3r-1)n\log (2)}{-(-4rn\log (2)+4p+2)}=1-\frac{3r-1}{4r}=\frac{r+1}{4r}$.

Clearly, the measurable entropy  $h_{\mu} (\sigma )$ on $X$, is equal to
$h_\mu(\sigma_1)+h_\mu (\sigma_2)=2\ln (2)$
and $h_\mu (F)=h_\mu (\sigma_1 )+h_\mu (\sigma_2^r)=(r+1)\ln (2)$.

 Finally  we can check the equality given in the previous subsection:
$$
h_{\mu} (F_e)=h_{\mu} (\sigma )\times 2r\times  M^*(r) =2\log (2)\times 2r\times \frac{r+1}{4r}=(r+1)\log (2).
$$


Remark that using Theorem 1, we obtain the strict following inequality
$h_\mu (F_e)\le (I_\mu^++I_\mu^-)h_\mu (\sigma)=r\times 2\ln (2)=2r\ln (2)$ with
$I_\mu^+=0$ and $I_\mu^-=r$.

In the following, we are going to extends the equality for a shift ergodic measure.


\subsection{The equalities and inequalities for a shift ergodic measure}

To extend the results of the previous section, we use the property that each cylinders defined
by fixing the same number of coordinates have a similar weight for a shift-ergodic measure.
This is due to the Shannon-McMillan-Breiman Theorem used for the shift action.

Denote by $\alpha_p$ the partition of $A^\ZZ$ into cylinders $[x(-p ,p)]_{-p}$ where
$x$ any point in $A^\ZZ$ and by
$^{\{n\}}\alpha_{p}^{F}$ the partition $\alpha_p\vee F^{-1}\alpha_p\ldots F^{-n}\alpha_p$
where $F :A^\ZZ \to A^\ZZ$ is a cellular automaton.
Then call $^{\{n\}}\alpha_p^{F}(x)$ the element of the partition $^{\{n\}}\alpha_p^{F}$ that contains
the point $x$.
Recall that for a two-sided subshift, the shift $\sigma$, is a bijective map.
For the bijective maps $\sigma$, we write $^G\alpha_p^{\sigma}(x)$ the element of the partition
$\vee_{i\in G}\; \sigma^{i}\alpha_p$  which contains the point $x$.
The set $G$ is a finite interval of $\ZZ$.
It follows  that each set $^G\alpha_p^{\sigma}(x)$ are cylinders set, while the set $^{\{n\}}\alpha_p^{F}(x)$ are
finite union of cylinders set.

Let $\mu$ be a shift  ergodic measure.
Given a real $0<\delta <1$ and two strictly increasing sequence maps $g^+_n$ and
 $g^-_n$ (from $A^\ZZ$ to $\NN$ with $\lim_{n\to\infty}g_n^+(x)+g_n^-(x)=+\infty$ for
 $\mu$ almost all $x$), set  $G_n(x)=\{-g^-_n(x),..., -1, 0, ...g^+_n(x)\}$,
 $|G_n(x)|=g_n^-(x)+g_n^+(x)+1$ and define

$$
\eta_{n,\delta}
=\max   \left\{
\begin{array}{ll}
 &\epsilon\in \RR:\; \exists S\subset X\; \mu (S)\ge 1-\delta \mbox{ and $S$ satisfies }\cr
&\forall x\in S, i\ge n:\;  \left\vert\frac{-\log \mu \left(^{G_n(x)}\alpha_p^{\sigma}(x)\right)}{|G_n(x)|}
-h_\mu (\sigma )
 \right\vert <\epsilon
\end{array}
\right\}
$$
and for any $n\in\NN$
$$
X^{G}_{n,\delta ,p}=\left\{ x\in A^\ZZ :\; \left\vert\frac{-\log \mu \left(^{G_n(x)}\alpha_p^{\sigma}(x)\right)}
{|G_n(x)|}-h_\mu (\sigma ) \right\vert \le\eta_{n,\delta} \right\}.
$$

Remark that from the Shannon-Breiman-McMillan Theorem, for all $0\le\delta\le 1$ we obtain that
 $\lim_{n\to\infty}\eta_{n,\delta}=0$ and $ \lim_{n\to\infty}$ $\mu (X^{G}_{n,\delta ,p})$ $=1$
for all sequence of application $G$.

\begin{defi}
Let define

$$
\left \langle T^G_{n,p}(x)\right\rangle  =\left\{
\begin{array}{l}
w\in A^{G_n(x)}\mbox{ such that } \exists y\in ^{\{n\}}\alpha_p^\sigma (x)\cr
 \mbox{ which verify }
y(-g_n^-(x),g_n^+(x))=w
\end{array}\right\}
$$
and
 $$
T^G_{n,p}(x) =\{y\in A^\ZZ, \mbox{ such that }y(-p,p)\in \langle T_n^G(x)\rangle \}
$$

Let  $T^G_{n,\delta ,p}(x)$ be the finite union of cylinders $[u]_{-|u|}\in A^\ZZ$,
with  $u\in A^{|G_n(x)|}$ such that
$$
T^{G}_{n,\delta ,p}(x)=T^G_{n,p}(x)\cap X^{G}_{n,\delta ,p}
$$
and
$$
\left \langle T^G_{n,\delta ,p}(x)\right\rangle  =\left\{
\begin{array}{l}
w=y(-g_n^-(x),g_n^+(x)) \mbox{ with } y\in T^{G}_{n,\delta ,p}(x)
\end{array}\right\}.
$$

\end{defi}

\begin{rem}
By taking the intersection with $X^{G}_{n,\delta}$, we keep the "good" cylinders which
have almost all the same weight.
Recall (see section 3.1) that for the uniform measure $\mu_u$  we have
$$
\mu_u\left(^{\{n\}}\alpha_p^{F}(x)\right)=\# T^{rn}_{n,p}(x)\times
\mu_u \left(^{\pm rn}\alpha_p^{\sigma}(x)\right)=\# T^{rn}_{n,p}(x)\times A^{-2rn-1}.
$$
In this case we can use $T^G_{n,p}(x)$ instead of $T^G_{n,\delta ,p}(x)$.
\end{rem}

Let $I^{+,*}_n(x)=\max_{y\in \alpha_p^F(x)}I^+_n(y)$ and
 $I^{-,*}_n(x)=\max_{y\in \alpha_p^F(x)}I^-_n(y)$.

Recall  that $g^+_n$ and $g^-_n$ are two increasing sequences of maps from $X$ to $\NN$
whom limit is $+\infty$. These maps are the base of the definition of the sequence of applications $G$.
\begin{lem}
If $g^+_n\ge I^{+*}_n$ and $g^-_n\ge I^{-*}_n$ then
$$
\mu (^{\{n\}}\alpha_p^{F}(x))=\sum_{y\in ^{\{n\}}\alpha_p^{F}(x)}
\mu (^{G_n(x)}\alpha_p^{\sigma}(y)).
$$
In this case we write that $G_n$ satify the condition (*).
\end{lem}
{\it Proof}


From  \cite{ti2000}  (Proposition 5.1) one has
 $$
^{(I_n^+(x)+p,I_n^-(x)+p)}\alpha_p^\sigma(x)\subset ^{\{n\}}\alpha_p^F(x).
$$
Since $g_n^+\ge I_n^{+,*}$ and $g_n^-\ge I_n^{-,*}$, then each cylinder set
 $^{G_n(x)}\alpha_p^\sigma (y)$  is a subset of $^{\{n\}}\alpha_p^F(x)$
  when  $y\in ^{\{n\}}\alpha_p^F(x)$.

Therefore $\{^{G_n(x)}\alpha_p^\sigma (y)\mbox{ with }y\in ^{\{n\}}\alpha_p^F(x)\}$ is a partition of
$^{\{n\}}\alpha_p^{F}(x)$.
\hfill$\Box$

\medskip

More generaly, remark that  $\mu (^{\{n\}}\alpha_p^{F}(x))\le \sum_{y\in ^{\{n\}}\alpha_p^{F}(x)}
\mu (^{G_n(x)}\alpha_p^{\sigma}(y))$ and $I_n^{\pm,*}(x)\le rn $ for all $x\in A^\ZZ$.

\begin{lem}
If $\mu$ is a shift-ergodic measure, for $\mu$-almost all points $x$  one has
\small
\begin{equation}
\lim_{n\to\infty}\frac{1}{n}\log \mu (^{\{n\}}\alpha_p^{F}(x))
\le\limsup_{\delta\to 0}\limsup_{n\to\infty}\frac{1}{n}\log \left(\mu (^{G_n(x)}\alpha_p^{\sigma}(x))
\times \# \langle T^{G}_{n,\delta ,p}(x)\rangle \right).
\end{equation}
\normalsize
If $G_n$ verify (*) the inequality becomes an equality and we obtain the convergence
with respect to the  variables $n$ an $\delta$.
\end{lem}
{\it Proof}

For all point $x$,  we have $\mu (^{\{n\}}\alpha_p^{F}(x))\le\sum_{y\in ^{\{n\}}\alpha_p^{F}(x)}
\mu (^{G_n(x)}\alpha_p^{\sigma}(y))$ and from Lemma 1, the inequality becomes an equality if condition
(*) hold. Hence,  it remains to prove that for almost all points $x$ we have
\small
$$
\lim_{\delta\to 0}\lim_{n\to\infty}\left[
\begin{array}{l}
\frac{1}{n}\log \left( \mu (^{G_n(x)}\alpha_p^{\sigma}(x))
\times \# \langle T^{G}_{n,\delta ,p}(x)\rangle \right)\cr
-\frac{1}{n}\log\left( \sum_{y\in ^{\{n\}}\alpha_p^{F}(x)}
\mu (^{G_n(x)}\alpha_p^{\sigma}(y)) \right)
\end{array}
\right]=0.
$$
\normalsize

We will show that for any $0<\delta' <1$, there exists a  set $S\subset A^\ZZ$ of measure $1-\delta$ such that
if  $x\in S$, the sequence
\small
$$
\left(
\frac{1}{n}\log \left( \mu (^{G_n(x)}\alpha_p^{\sigma}(x))
\times \# \langle T^{G}_{n,\delta ,p}(x)\rangle \right)
-\frac{1}{n}\log\left( \sum_{y\in ^{\{n\}}\alpha_p^{F}(x)}
\mu (^{G_n(x)}\alpha_p^{\sigma}(y)) \right)\right)
$$
\normalsize
converge to 0 when $n$ goes to infinity.
Then we claim  that if $x$ belong to a set of measure $1-\delta '$ ($0<\delta ''<1$) we obtain
\small
\begin{equation}
\lim_{n\to\infty}\frac{1}{n}\left( \log \sum_{y\in ^{\{n\}}\alpha_p^{F}(x)}
\mu \left(^{G_n(x)}\alpha_p^{\sigma}(y)\right)-\log \sum_{y\in T^{G}_{n,\delta}(x)}
\mu \left(^{G_n(x)}\alpha_p^{\sigma}(y)\right)  \right) =0 .
\end{equation}
\normalsize
Let's prove this claim.
Fix $k>1$ and for all integer $n>0$, denote by  $Y_{n,\delta}$
  the set of points $x$ such that
$$
\frac{\sum_{y\in ^{\{n\}}\alpha_p^{F}(x)}
\mu \left(^{G_n(x)}\alpha_p^{\sigma}(y)\right)}
{\sum_{y\in T^{G}_{n,\delta}(x)}
\mu \left(^{G_n(x)}\alpha_p^{\sigma}(y)\right)}\le k.
$$

It is clear that all  points in $Y_\delta=\lim_{n\to\infty}\cap_{i=0}^n \cup_{j=i}^\infty Y_{j,\delta}$
 verify  equality (2). Hence, in order to prove the  claim,  we need to show that for all $n\in\NN$ we have
$\mu (^\complement Y_{n,\delta})\le \delta ''$ where $^\complement Y_{n,\delta}$ is the complementary set of
 $Y_{n,\delta}$.

The following  complete the proof of the claim:
Since
$$
^\complement Y_{n,\delta}\subset
 \left( ^{G_n(x)}\alpha_p^\sigma (x)-T_{n,\delta ,p}^{G}(x), \, x\in ^\complement Y_{n,\delta}\right)
\subset  \left( A^\ZZ-X^G_{n,\delta ,p} \right),
$$
then
$$
\frac{k}{k+1}\mu \left( ^{G_n(x)}\alpha_p^\sigma (x), \, x\in ^\complement Y_{n,\delta}\right)
\le \mu \left( A^\ZZ-X^G_{n,\delta ,p} \right),
$$
which implies that
$$
\frac{k}{k+1}\mu (^\complement Y_{n,\delta})\le \left( A^\ZZ-X^G_{n,\delta ,p} \right)
\mbox{ and } \mu (^\complement Y_{n,\delta})\le \frac{k+1}{k}\delta:=\delta'' .
$$

Remark that if $\mu ( Y_{n,\delta})\le 1-\delta ''$ for all integer $n$, then
$\mu ( Y_{\delta})\le 1-\delta ''$.

 Next, it remains  to prove that for all fixed $0<\delta <1$ and point  $x$ in
 $X^G_{\delta ,p}=\lim_{n\to\infty}\cap_{i=0}^n \cup_{j=i}^\infty X^{G}_{j,\delta ,p}$,
verify
\begin{equation}
\lim_{n\to\infty}\frac{1}{n}\left(
\log (\hskip -.2 cm \sum_{y\in T^G_{n,\delta}(x)} \mu (^{G_n(x)}\alpha_p^{\sigma}(y))
-\log \left( \#\langle T^G_{n,\delta ,p}(x)\rangle\times \mu (^{G_n(x)}\alpha_p^{\sigma}(x)\right)\right)=0.
\end{equation}

Remark that  if  $x\in X^{G}_{n,\delta ,p}$ and $y\in T^{G}_{n,\delta ,p}(x)$ one has
$
\left\vert \frac{\mu (^{G_n(x)}\alpha_p^{\sigma}(y))}{\mu (^{G_n(x)}\alpha_p^{\sigma}(x))}
\right\vert\le e^{n\eta_{n,\delta}}.
$
It follows that

 $$
\frac{1}{n}\left(
\log ( \sum_{y\in T^G_{n,\delta}(x)} \mu (^{G_n(x)}\alpha_p^{\sigma}(y))
-\log \left( \#\langle T^G_{n,\delta ,p}(x)\rangle\times \mu (^{G_n(x)}\alpha_p^{\sigma}(x)\right)\right)=\eta_{n,\delta}.
$$

Since for $0<\delta<1$ the  Shannon-Breiman-McMillan  Theorem  gives that
$\lim_{n\to\infty}\eta_{n,\delta}=0$, so we have proved the equation (3).

Arguing that equation (1) is true for all $x\in X_\delta^G\cap Y_\delta=:S$ and
taking into consideration that $\mu (X_\delta^G\cap Y_\delta)\ge 1-\delta-\frac{k+1}{k}\delta=:\delta'$,
we can conclude letting $\delta\to 0$.

From the Shannon-Breiman-McMillan Theorem,
$\left(\frac{1}{n}\log \mu (^{\{n\}}\alpha_p^{F}(x))\right)_{n\in\NN}$ converge almost everywhere when
$\mu$ is an invariant measure, this implies that when $G_n$ verify condition (*) we obtain
$$
\lim_{n\to\infty}\frac{1}{n}\log \mu (^{\{n\}}\alpha_p^{F}(x))
=\lim_{\delta\to 0}\lim_{n\to\infty}\frac{1}{n}\log \left(\mu (^{G_n(x)}\alpha_p^{\sigma}(x))
\times \# \langle T^{G}_{n,\delta ,p}(x)\rangle \right).
$$

\hfill$\Box$

\begin{defi}

For each integers $p>0$, $n>0$, real $0<\delta <1$ and double sequences of functions
 $G=((g_n^-(x)),(g_n^+(x)))$ from $A^\ZZ$ to $\NN$ we write
$$
M_{n,p,\delta}(G)=\int_{A^\ZZ} 1-\frac{\log
\# \langle T^{G_n(x)}_{n,\delta ,p}(x)\rangle}{-\log \mu (^{G_n(x)}\alpha_p^{\sigma}(x))}d\mu (x)
$$

and
$$
M_{\mu} (G)=\sup_{p}\limsup_{\delta\to 0}\limsup_{n\to\infty}M_{n,p,\delta}(G).
$$
We call  $M_\mu (G)$  the average flow of information at speed $G$.

We denote by $I^{+,*}$ the value $\sup_{A^\ZZ}\limsup_{n\to\infty}\frac{I_n^{+,*}(x)}{n}$
and by $I^{-,*}$ the value $
\sup_{A^\ZZ}\limsup_{n\to\infty}\frac{I_n^{-,*}(x)}{n}$.

\end{defi}


\begin{rem}
In the following theorem
$
M_{\mu} (G)=\sup_{p}\lim_{\delta\to 0}\lim_{n\to\infty}M_{n,p,\delta}(G)
$
 when the velocity $v$ is greater than $(I^{+,*},I^{-,*})$.
Remark that $I^{+,*}, I^{-,*}\le r$ where $r$ is the radius of the cellular automaton $F$ we consider.
\end{rem}

\begin{thm}
If  $\mu$ is a $\sigma$-ergodic and $F$-invariant measure and
$v=\left((\lceil v^-n\rceil )\right.$ $\left.,(\lceil v^+n\rceil)\right)$ a double sequence of integers
we have the following properties :
\medskip

(i) $h_\mu (F)=h_\mu (\sigma )\times (v^++v^-)\times M_\mu (v)$, if
$v^+\ge I^{+,*}$ and  $v^- \ge I^{-,*}$.

\medskip

(ii) $h_\mu (F)\ge h_\mu (\sigma )\times (v^++v^-)\times M_\mu (v)$.
\medskip


(iii)
$
h_\mu (F)=h_\mu (\sigma )\times M
$
where
$$
M=sup_p\int_{A^\ZZ}\lim_{\delta\to 0} \lim_{n\to\infty}
M_{n,\delta,p}(I^{+*}_n(x),I_n^{-*}(x))\times \frac{I^{\pm*}_n(x)}{n} d\mu (x).
$$
and $I^{\pm*}_n(x)=I^{+,*}_n(x)+I_n^{-,*}(x)$.
\normalsize
\end{thm}

{\it Proof}


Applying the Shannon-McMillan-Breiman Theorem (probabilistic version) to $F$ we obtain
$$
h_\mu (F,\alpha_p )=\int  \lim_{n\to\infty}-\frac{1}{n}\log \mu (^{\{n\}}\alpha_p^{F}(x))d\mu (x)
$$

In order to simplify the first part of the proof we suppose condition
(*) (which is the case in part (i) and (iii)). In the more genaral case, we have
limits superior  instead of limits and inequalities instead of equalities.

From Lemma 2, for almost all $x$ we have
\small
$$
\lim_{n\to\infty}\frac{1}{n}\log \mu (^{\{n\}}\alpha_p^{F}(x))=
\lim_{\delta\to 0}\lim_{n\to\infty}\frac{1}{n}\log \left(\mu \left(^{G_n(x)}\alpha_p^{\sigma}(x)\right)\times
\log (\# \langle T^{G_n(x)}_{n,\delta ,p}(x)\rangle \right)
$$
\normalsize
$$
=
 \lim_{\delta\to 0}\lim_{n\to\infty}\frac{1}{n}\left(\log \mu (^{G_n(x)}\alpha_p^{\sigma}(x))+
\log (\# \langle T^{G_n(x)}_{n,\delta ,p}(x)\rangle )\right).
$$
We can rewrite the inequality as

$\lim_{n\to\infty}-\frac{1}{n}\log \mu (^{\{n\}}\alpha_p^{F}(x))=$
$$
\lim_{\delta\to 0}\lim_{n\to\infty}|G_n(x)|\frac{-\log \mu (^{G_n(x)}\alpha_p^{\sigma}(x))}
{|G_n(x)|n}- \frac{\log (\# \langle T^{G_n(x)}_{n,\delta ,p}(x))\rangle }{n}
$$

and we obtain that $\lim_{n\to\infty}\frac{1}{n}\log \mu (^{\{n\}}\alpha_p^{F}(x))=$
$$
\lim_{\delta\to 0}\lim_{n\to\infty} \frac{|G_n(x)|}{n}\times \frac{-\log \mu
 (^{G_n(x)}\alpha_p^{\sigma}(x))}{|G_n(x)|}\times
\left (1-\frac{\log (\# \langle T^{G_n(x)\rangle }_{n,\delta ,p}(x))}
{-\log \mu (^{G_n(x)}\alpha_p^{\sigma}(x))}\right).
$$

Hence from the dominated convergence theorem one has

$$
\hskip -11 true cm h_\mu (F,\alpha_p)=
$$
\small
$$
 \lim_{\delta\to 0}\lim_{n\to\infty}\int  \frac{-\log \mu
 (^{G_n(x)}\alpha_p^{\sigma}(x))}{|G_n(x)|}\times
\frac{|G_n(x)|}{n}\times
\left (1-\frac{\log \left(\# \langle T^{G_n(x)}_{n\delta ,p}(x)\rangle\right) }
{-\log \mu (^{G_n(x)}\alpha_p^{\sigma}(x))}\right)
d\mu (x)
$$
\normalsize

Using the extented version of Shannon-Breiman-McMillan Theorem for $\sigma$ (see \cite{Or83}), we obtain
that
$h_\mu (F,\alpha_p)$
$$=h_\mu (\sigma ,\alpha_p )\times \lim_{\delta\to 0}\lim_{n\to\infty}\int
\frac{|G_n(x)|}{n}\times
\left (1-\frac{\log \left(\# \langle T^{G_n(x)}_{n\delta ,p}(x)\rangle\right) }
{-\log \mu (^{G_n(x)}\alpha_p^{\sigma}(x))}\right)
d\mu (x).
$$


Proof of (i)

If $G(n)=(\lceil v^-n\rceil ,\lceil v^+n\rceil)$ and $v^+\ge I^{+,*}$, $v^-\ge I^{-,*}$ then
 the equality in Lemma 2 becomes an equality, the limsup in $n$ is a limit and we can write that

$h_\mu (F,\alpha_p)$
$$
=h_\mu (\sigma ,\alpha_p)\times (v^++v^-)\times
\lim_{\delta\to 0}\lim_{n\to\infty}
\int 1-\frac{\log \left(\# \langle T^{G}_{n\delta ,p}(x)\rangle\right) }
{-\log \mu (^{(\lceil v^-n\rceil ,\lceil v^+n\rceil)}\alpha_p^{\sigma}(x))} d\mu (x).
$$

and finally
$$
h_\mu (F)=h_\mu (\sigma )\times (v^++v^-)\times M_\mu (v)
$$

where
$$
M_\mu (v)=\sup_{p}\lim_{\delta\to 0}\lim_{n\to\infty}\int_{A^\ZZ} 1-\frac{\log
\# \langle T^{G_n(x)}_{n,\delta ,p}(x)\rangle}{-\log \mu (^{G_n(x)}\alpha_p^{\sigma}(x))}d\mu (x).
$$
Remark that the Shannon-Breiman-McMillan theorem implies that the density flow can be defined
as a double limit in that case.
\hfill$\Box$

Proof of (ii)

Same proof but using Lemma 2 without condition (*).
Remark that, in that case the density flow is defined thanks to  limits superior instead of simple limits.

\hfill$\Box$

Proof of (iii)

From \cite{ti2000} if $F$ has equicontinuous points then $\frac{I^+_n(x)}{n}$ and $\frac{I^+_n(x)}{n}$
are almost everywhere  bounded and $h_\mu (F)=0$.

When  $F$ is sensitive (no equicontinuous points) then for all $x$ we have
 $\liminf_{n\to\infty}\frac{I_n^+(x)}{n}+\frac{I_n^-(x)}{n}=+\infty$.
Clearly we have the same properties for $I_n^{+,*}(x)$ and $I_n^{-,*}(x)$.
Using Lemma 2 with condition (*) we establish that

$$
\hskip -8 cm h_\mu (F,\alpha_p)
=h_\mu (\sigma ,\alpha_p )\times
$$
$$
 \lim_{\delta\to 0}\lim_{n\to\infty}\int
 \frac{I^{+*}_n(x)+I^{-*}_n(x)}{n}\times
\left (1-\frac{\log \left(\# \langle T^{I^{*}}_{n,\delta ,p}(x))\rangle\right) }
{-\log \mu (^{I^*_n (x)}\alpha_p^{\sigma}(x))}\right)
d\mu (x)
$$
where $I^{*}_n(x)=(I_n^{+,*}(x),I_n^{-,*}(x))$.

\hfill$\Box$


\begin{ques}
Obviously for all $x$,  one has $\liminf_{n\to\infty}\frac{I_n^{\pm,*}(x)}{n}$ is greater than
$\liminf_{n\to\infty}\frac{I_n^{\pm}(x)}{n}$.
What are the other relations between these two functions.
\end{ques}

\begin{rem}
\mbox{  }

$\bullet$
Since  $\liminf_{n\to\infty}\frac{I_n^{\pm,*}(x)}{n}$ could be
equal to zero almost everywhere, the part (iii) of Theorem 2
 gives  sens to the study of the density flow for sublinear speed.

$\bullet$ If the constants $I_\mu^+$ and $I_\mu^-$ represent
the average speed of the faster perturbations, in the point of view of the density flow,
the ratio $\frac{h_\mu (F)}{h_\mu (\sigma )}=(v^++v^-)M_\mu (v)$
gives the sum of the
average  over all the points of the averages perturbations's speeds when $v^\pm\ge r$.
\end{rem}



\subsection{Some Examples}

\textsc{Example 1}

We use the same automaton $F_e$ than in the
subsection 3.1.2 (based on the shifts $\sigma_1$ on $X_1=\{0,1\}^\ZZ$ and $\sigma_2$ on
$X_2=\{0,1\}^\ZZ$)
 but we consider the non uniform Bernouilli measure
 $\mu_B=\mu_{B1}\times\mu_{B2}=B(\frac{1}{3},\frac{2}{3})\times B(\frac{1}{3},\frac{2}{3}) $.
The entropy on the two shifts
 $h_{\mu_{B1}} (\sigma_1)=h_{\mu_{B2}} (\sigma_2)=\frac{1}{3}\log (3)+\frac{2}{3}\log (\frac{3}{2})$
 and $h_{\mu_B}(\sigma )=h_\mu (\sigma_1)+h_\mu (\sigma_2)=
2(\frac{1}{3}\log (3)+\frac{2}{3}\log (\frac{3}{2}))$.
Considering a set of words $\langle T^{rn}_{n,p}(x)\rangle$ instead of set of cylinders $T^{rn}_{n,p}(x)$,
 and using the same arguments than in section 3.1.2,
 we obtain that
 $\#\langle T^{rn}_{n,p}(x)\rangle =4^{rn+2p+1}\times 2^{(r-1)n}=2^{(3r-1)n}$.

For $\mu_B$ almost all the point $x$, the number $T^{rn}_{n,\delta ,p}(x)$ of "good cylinders"
(in $X^{rn}_{n, \delta ,p}$)
 among the $2^{(3r-1)n}=\# \langle T^{rn}_{n,p}(x)\rangle$ possible words is approximatively
 (when $n$ large enough):
 $\# \langle T^{rn}_{n,\delta ,p}\rangle\approx e^{h_\mu (\sigma )(3r-1)n}$ (a consequence of Shannon-Breiman-McMillan
 Theorem, see \cite{WA} page 95)
and for almost all $x$, $-\log \mu (^{rn}\alpha_p^\sigma (x))\approx e^{(4rn+4p+2)h_\mu (\sigma_1)}$.

Therefore we get,
$$
M_\mu (r)=\lim_{\delta\to 0}\lim_{n\to\infty}
1-\frac{\log \# T^{rn}_{n,\delta ,p}(x)}{-\log \mu (^{rn}\alpha_p^\sigma (x))},
$$

$$
M_\mu (r)=1-\frac{3r-1}{4r}=\frac{r+1}{4r}.
$$

Finaly we obtain
$$
h_{\mu_B}(F)=M_\mu (r)\times h_{\mu_B}(\sigma )\times 2r=\frac{r+1}{4r}
\times 2r\times 2(\frac{1}{3}\log (3)+\frac{2}{3}\log (\frac{3}{2})),
$$
$$
=(r+1)(\frac{1}{3}\log (3)+\frac{2}{3}\log (\frac{3}{2}))=h_{\mu_B} (\sigma_1)+
h_{\mu_B}(\sigma_2^r).
$$
\bigskip

\textsc{Examples 2}


Let $X_2$ be a Sturmien minimal subset of $\{0,1\}^\ZZ$. Let $\mu^s_2$ be the unique ergodic
invariant measure on
$(X_2,\sigma_2 )$ where $\sigma_2$ is the shift on $X_2$. We call  $\mu_1$,
the uniform measure on $X_1=\{0,1\}^\ZZ$ and we denote by $x=x^1\times x^2$ any point in $X=X_1\times X_2$.
Consider  the cellular automaton
$F =Id_1\times\sigma_2^r$  where $Id_1$ is the identity map on $X_1$, $r$ any positive integer and
denote by $\mu$ the measure $\mu_1\times\mu_2^s$ defined on $X$.
Using well known results,
we have  $h_{\mu^s_2}(\sigma_2)=0$, $h_{\mu_1}(Id_1)=0$ and consequently $h_\mu (F)=0$.
From theorem 2, we have
 $M_\mu (r)=0$.
Clearly the Lyapunov exponents $I_\mu^++I_\mu^-=I^+_\mu=r$, so in this case there exist perturbations moving at a
positive speed but the average flow is equal to 0.

If we modify slightly the last example and we replace the identity map on $X_1$ by  $r$ iterations of
the shift $\sigma_1$ on
$X_1$ and the shift $\sigma_2^r$ by the identity, we obtain a couple
(cellular automaton $F_2=\sigma_1^r\times Id_2$, measure $\mu$)
 such that the inequality of Theorem 1 becomes an equality :
  $h_\mu (F_2)=h_\mu (\sigma )\times (I^+_\mu +I^-_\mu)
=r\log (2)$ and in this case $M_\mu (r,0)=1$ (the maximum value of $M_\mu (v)$).

\subsection{The density flow and some subclasses of cellular automata}

Clearly, in some cases like the shifts or more generaly bipermutative cellular automata
(see \cite{GI88}), all the perturbations
 move at  speed $r$ which implies that $M_\mu(r)=1$ for all shift ergodic and $F$-invariant measure $\mu$.
 Sometimes, like in the last example $M_\mu(r)=1$ and almost all perturbation move at speed $r$.
The next Proposition shows that  there exists  a kind of ``minimum speed'' for all the
perturbations in the case of the positively expansive CA.

\begin{pro}
If $F$ is a positively expansive CA and $\mu$ a shift ergodic and $F$-invariant measure,
 then there exist a real $s>0$ such that
$M_\mu (s)=1$
\end{pro}
{\it Proof}

From \cite{BT2006}(Proposition 2), all positively expansive CA has positive pointwise Lyapunov exponents
for all point in $A^\ZZ$.
More precisely there exist a  real $\eta >0$ such that
 $\liminf_{n\to\infty}\frac{I_n^+(x)}{n}>\eta $ and $\liminf_{n\to\infty}\frac{I_n^-(x)}{n}> \eta$
for all $x\in A^\ZZ$.

If the velocity of the perturbations  is $(\eta, \eta)$, there would exist an integer
 $N>0$ such that for all $n\ge N$ we have $\frac{I_n^+(x)}{n}+\frac{I_n^-(x)}{n}> 2\eta$.

In this case for all integer $n\ge N$, $p>0$,  real $0<\delta <1$ and $x\in A^\ZZ$, we obtain
$\#\langle T^{(\eta n,\eta n)}_{n, \delta ,p}(x)\rangle =1$
which implies that  $\log \left (\# T^{(\eta n, \eta n)}_{n, \delta ,p}(x)\right) =0$ and
consequently $M_\mu ((\eta , \eta ))=1$.

\hfill$\Box$

\begin{ques}
From Theorem 2, the function $M_\mu (s)$ from $[r,+\infty ) \to [0,1]$ is clearly continuous and decreasing
but what happens in the interval $[0,1]$?
\end{ques}








Since $\mu$-expansiveness property is a measurable equivalent to positive expansiveness, we can
wonder what values can take the density flow in that case.
We will see that there is no such strong relation with the density flow but for the
complementary class of CA with $\mu$-equicontinuous points, the density flow is always equal to zero.

From $\cite{GI87}$ the existence of  equicontinuous points is equivalent to the existence of
blocking words that stop completly the propagation of the perturbations. From \cite{tippt2006},  there
exist an example of couple (cellular automaton, invariant and shift ergodic measure $\mu$)
 without equicontinuous points such that there exists $\mu$-equicontinuous points for the cellular automaton
 action restricted to the topological
support of $\mu$. This example shows clearly than the dynamic of CA with $\mu$-equicontinuous points is clearly
more general than those with equicontinuous points when $\mu$ is an invariant measure.
The next Proposition shows that the density of flows of the perturbations
that can move or not move to infinity (strictly $\mu$-equicontinuous or equicontinuous case)
is equal to zero taking in consideration a linear or a sublinear speed.

\begin{pro}
If $F$ is a $\mu$-equicontinuous CA, then for all couple of sequence $G_n=(g_n^+,g_n^-)$ such that
$\lim_{n\to\infty}g_n^++g_n^-=+\infty$ we have $M_\mu(G)=0$.
\end{pro}
{\it Proof}

Suppose there exist a couple of sequence $G_n=(g_n^+,g_n^-)$ such that $M_\mu(G)$ $>0$ with
$\lim_{n\to\infty}g_n^++g_n^-=+\infty$.
It follows that there if we fixe $\delta$ and $p$, there exist a set $S$ of positive measure and a real
$0<\alpha <1$ such that
for all point $x\in S$ we have
$$
\lim_{n\to\infty}\log \# \langle T^{G}_{i(n),\delta ,p}(x)\rangle
<\alpha \lim_{n\to\infty}\left[ -\log \mu \left( ^{G_{i(n)}} \alpha_p^\sigma (x)\right)\right],
$$
where $i(n)$ is a subsequence that allows the convergence.

Then,  for $n$ large enough, there exist $\alpha '$ such that
 $$
 \log \# \langle T^{G}_{i(n),\delta ,p}(x)\rangle
<\alpha \left[ -\log \mu \left( ^{G_{i(n)}} \alpha_p^\sigma (x)\right)\right].
$$

Using the same argument as in the  proof of Lemma 2, we obtain that
for $n$ large enough, there exist a subset of $S$: $S'$ of positive measure such that
there exist a positive $k>1$ such that

$$
\frac{\sum_{y\in ^{\{n\}}\alpha_p^{F}(x)}
\mu \left(^{G_n}\alpha_p^{\sigma}(y)\right)}
{\sum_{y\in T^{G}_{n,\delta}(x)}
\mu \left(^{G_n}\alpha_p^{\sigma}(y)\right)}\le k.
$$

Since $\mu (^{\{n\}}\alpha_p^{F}(x))\le \sum_{y\in ^{\{n\}}\alpha_p^{F}(x)}
\mu (^{G_n}\alpha_p^{\sigma}(y))$ and when $n$ is large enough
$\mu \left( ^{G_n} \alpha_p^\sigma (x)\right)\approx$ $e^{|G_n|h_\mu (\sigma )}$, then
$\mu (^{\{n\}}\alpha_p^{F}(x))\le
(k+1)\sum_{y\in T^{G}_{n,\delta}(x)}
\mu \left(^{G_n}\alpha_p^{\sigma}(y)\right)
$.
It follows that
\small
$$
\mu (^{\{n\}}\alpha_p^{F}(x))\le (k+1)\# \langle T^{G}_{n,\delta ,p}(x)\rangle \times
\mu \left( ^{G_n} \alpha_p^\sigma (x)\right)
$$
$$
\le
(k+1)\mu \left( ^{G_n} \alpha_p^\sigma (x)\right)^{-\alpha}\mu \left( ^{G_n} \alpha_p^\sigma (x)\right).
$$
\normalsize
This implies that for a set of positive measure $\lim_{n\to\infty}\mu (^{\{n\}}\alpha_p^{F}(x))=0$
 $=\mu ( B_p(x))$. From \cite{GI87} and \cite{tippt2006}, if there exists a $\mu$-equicontinuous point,
 then there exist a set of measure 1 of these kind of points.
 Hence  there is no $\mu$-equicontinuous point and we can conclude taking the reverse asumption.

 \hfill$\Box$
\medskip

Using Proposition 3 (only for linear speed) and Theorem 2, it follows that

\begin{cor}
If $\mu$ is a shift ergodic and $F$-invariant measure, with $F$ a cellular automaton with
$\mu$-equicontinuous points, then the measurable entropy $h_\mu (F)=0$.
\end{cor}












\begin{rem}
If $F$ is a $\mu$-expansive CA then clearly for $p$,  $\delta$ fixed and  all $x$ in
a set of positive measure we have
$$
\lim_{n\to\infty}\# T_{n,p,\delta}^{I^{*}}(x)\times \mu (^{I^{*}_n(x)}\alpha_p^\sigma (x))=0
$$ where
$I^*=(I^{+,*}_n(x),I=(I^{-,*}_n(x))$. Remark that this condition is weaker to the property $M_\mu (I^*)>0$.
\end{rem}

\begin{ques}
Is it possible that a dynamical system  (cellular automaton $F$, invariant measure $\mu$),
 verify that $M_\mu (I^*)>0$ and $I^{+,*}_\mu +I^{-,*}_\mu =0$?
 The sum of the exponents $I^{+,*}_\mu$ and $I^{-,*}_\mu$ is defined by
$$
I^{+,*}_\mu +I^{-,*}_\mu=\limsup_{n\to\infty}\int \frac{I^{+,*}_n(x)+I^{-,*}(x)}{n}d\mu (x)=0.
$$
From Proposition 2, this CA have to be $\mu$-expansive.
In \cite{BT2006}, there is an example of $\mu$-expansive CA with $I^{+}_\mu +I^{-}_\mu=0$
but $M_\mu (I^*)$ and $I^{+,*}_\mu$ and $I^{-,*}_\mu$ are not known.
\end{ques}

We finish with some more general questions.

\begin{ques}
For $h_\mu (F,\alpha_p)$, the bound given by Theorem 1 using the Lyapunov exponents $I_\mu^++I_\mu^-$
did not depends on the choice of the finite partition $\alpha_p$.
In what cases the value of $M_\mu (v)$ does not depends on a supremum value over all the partition $\alpha_p$?
\end{ques}





\begin{ques}
For each cellular automaton $F$, it is possible to defined the function $M_\mu (v)$ when $\mu$
is shift ergodic and non $F$-invariant. What is the meaning of $M_\mu (v)$ in this case?
\end{ques}

\begin{ques}
What kind of results can be done in the topological case?
\end{ques}


\end{document}